\documentclass[final]{siamltex}

\usepackage[parfill]{parskip}    
\usepackage{amsmath,amssymb,latexsym,enumerate,mathrsfs}
\usepackage{hyperref}
\usepackage{array}
\usepackage{times}
\usepackage[]{mdframed}
\usepackage{epstopdf}
\usepackage{subfig}
\usepackage{graphicx}
\usepackage{hyperref} 
\usepackage{multicol}
\usepackage{framed}
\usepackage{moreverb}
\usepackage{marvosym}
\usepackage{color,soul}
\definecolor{lightblue}{rgb}{.90,.95,1}
\sethlcolor{yellow}

\DeclareMathOperator*{\argmax}{arg\,max}

\def\qed{\hfill $\Box$}

\title{On improving generalization in a class of learning problems with the method of small parameters for weakly-controlled optimal gradient systems}

\author{Getachew K. Befekadu}

\begin{document}
\maketitle

\renewcommand{\thefootnote}{\arabic{footnote}}

\begin{abstract}
In this paper, we provide a mathematical framework for improving generalization in a class of learning problems which is related to point estimations for modeling of high-dimensional nonlinear functions. In particular, we consider a variational problem for a weakly-controlled gradient system, whose control input enters into the system dynamics as a coefficient to a nonlinear term which is scaled by a small parameter. Here, the optimization problem consists of a cost functional, which is associated with how to gauge the quality of the estimated model parameters at a certain fixed final time w.r.t. the model validating dataset, while the weakly-controlled gradient system, whose the time-evolution is guided by the model training dataset and its perturbed version with small random noise. Using the perturbation theory, we provide results that will allow us to solve a sequence of optimization problems, i.e., a set of decomposed optimization problems, so as to aggregate the corresponding approximate optimal solutions that are reasonably sufficient for improving generalization in such a class of learning problems. Moreover, we also provide an estimate for the rate of convergence for such approximate optimal solutions. Finally, we present some numerical results for a typical case of nonlinear regression problem.
\end{abstract}


\section{Statement of the problem} \label{S1}
Consider the following optimal control problem (which is related to a weakly-controlled gradient optimal system with a small parameter)
\begin{align}
J^{\epsilon}[u] &= \Phi\bigl(\theta^{\epsilon}(T), \mathcal{Z}^{(2)}\bigr) \quad \to \quad \min_{u(t) \in U ~ \text{on} ~ t \in [0, T]} \label{Eq1.1}\\ 
             & \text{s.t.}  \notag\\
              \dot{\theta}^{\epsilon}(t) &= - \nabla J_{0} \bigl(\theta^{\epsilon}(t), \mathcal{Z}^{(1)}\bigr) + \epsilon u(t) B\bigl(\theta^{\epsilon}(t), \tilde{\mathcal{Z}}^{(1)}\bigr), \quad \theta^{\epsilon}(0)=\theta_0, \label{Eq1.2}
\end{align}
where the statement of the problem consists of the following core concepts and general assumptions:
\begin{enumerate} [(a).]
\item {\bf Datasets}: We are given two datasets, i.e., $\mathcal{Z}^{(k)} = \bigl\{ (x_i^{(k)}, y_i^{(k)})\bigr\}_{i=1}^{m_k}$, each with data size of $m_k$, for $k=1, 2$. These datasets, i.e., $\mathcal{Z}^{(1)}$ and $\mathcal{Z}^{(2)}$, may be generated from a given original dataset $\mathcal{Z}^{(0)} =\bigl\{ (x_i^{(0)}, y_i^{(0)})\bigr\}_{i=1}^{m_0}$ by means of bootstrapping with/without replacement. Here, we assume that the first dataset $\mathcal{Z}^{(1)}=\bigl\{ (x_i^{(1)}, y_i^{(1)})\bigr\}_{i=1}^{m_1}$ will be used for model training purpose, while the second dataset $\mathcal{Z}^{(2)} = \bigl\{ (x_i^{(2)}, y_i^{(2)})\bigr\}_{i=1}^{m_2}$ will be used for evaluating the quality of the estimated model parameter. Moreover, the dataset $\tilde{\mathcal{Z}}^{(1)}$ (which is associated with the nonlinear term $B$ in the system dynamics) is obtained by adding small random noise, i.e., $\tilde{\mathcal{Z}}^{(1)} = \bigl\{ (x_i^{(1)}, \tilde{y}_i^{(1)})\bigr\}_{i=1}^{m_1}$, with $\tilde{y}_i^{(1)} = y_i^{(1)} + \varepsilon_i$ and $\varepsilon_i \sim \mathcal{N}(0, \sigma^2)$ (with small variance $\sigma^2$).
\item {\bf Learning via weakly-controlled gradient systems with a small parameter}: We are tasked to find for a parameter $\theta \in \Theta$, from a finite-dimensional parameter space $\mathbb{R}^p$ (i.e., $\Theta \subset \mathbb{R}^p$), such that the function $h_{\theta}(x) \in \mathcal{H}$, i.e., from a given class of hypothesis function space $\mathcal{H}$, describes best the corresponding model training dataset as well as predicts well with reasonable expectation on a different model validating dataset. Here, the search for an optimal parameter $\theta^{\ast} \in \Theta \subset \mathbb{R}^p$ can be associated with the weakly controlled-gradient system of Equation~\eqref{Eq1.2}, whose {\it time-evolution} is guided by the model training dataset $\mathcal{Z}^{(1)}$ and its perturbed version $\tilde{\mathcal{Z}}^{(1)}$, i.e.,
\begin{align*}
 \dot{\theta}^{\epsilon}(t) = - \nabla J_{0} \bigl(\theta^{\epsilon}(t), \mathcal{Z}^{(1)}\bigr) + \epsilon u(t) B\bigl(\theta^{\epsilon}(t), \tilde{\mathcal{Z}}^{(1)}\bigr), 
 \end{align*}
 where $J_0\bigl(\theta, \mathcal{Z}^{(1)}\bigr) = \frac{1}{m_1} \sum\nolimits_{i=1}^{m_1} {\ell} \bigl(h_{\theta}(x_i^{(1)}), y_i^{(1)} \bigr)$, and $\ell$ is a suitable loss function that quantifies the lack-of-fit between the model and the datasets. Moreover, $u(t)$ is a real-valued admissible control function from a set $U$ in one-dimensional space that enters into the system dynamics as a coefficient to the nonlinear term $B$.\footnote{Note that the control function $u(t)$ is admissible if it is measurable and $u(t) \in U$ for all $t \in [0, T]$.} The parameter $\epsilon$ is a small positive number and the nonlinear term $B$ is given by
\begin{align*}
 B\bigl(\theta, \tilde{\mathcal{Z}}^{(1)}\bigr) = \left [\bigl(\partial J_0 \bigl(\theta, \tilde{\mathcal{Z}}^{(1)}\bigr)/ \partial \theta_1\bigr)^2, \bigl(\partial J_0 \bigl(\theta, \tilde{\mathcal{Z}}^{(1)}\bigr)/ \partial\theta_2\bigr)^2,\ldots, \bigl(\partial J_0 \bigl(\theta, \tilde{\mathcal{Z}}^{(1)}\bigr)/ \partial \theta_p\bigr)^2 \right]^T 
\end{align*}
Note that the small random noise $\varepsilon \sim \mathcal{N}(0, \sigma^2)$ (with small variance $\sigma^2$) in the dataset $\tilde{\mathcal{Z}}^{(1)}$ provides a {\it dithering effect}, i.e., causing some distortion to the model training dataset $\mathcal{Z}^{(1)}$ so that the control $u(t)$ will have more effect on the learning dynamics.\footnote{In Equation~\eqref{Eq1.2} above, the weakly controlled-gradient system can be expressed as follows:
\begin{align*}
 \dot{\theta}_i^{\epsilon}(t) = - \partial J_{0} \bigl(\theta^{\epsilon}(t), \mathcal{Z}^{(1)}\bigr)/ \partial \theta_i + \epsilon u(t) \bigl(\partial J_0 \bigl(\theta, \tilde{\mathcal{Z}}^{(1)}\bigr)/ \partial \theta_i\bigr)^2, \quad i = 1,2, \ldots, p,
\end{align*}
where the control $u(t)$ enters into the system dynamics as a common coefficient to all nonlinear terms $\bigl(\partial J_0 \bigl(\theta, \tilde{\mathcal{Z}}^{(1)}\bigr)/ \partial \theta_i\bigr)^2, \tilde{\mathcal{Z}}^{(1)}\bigr)$, $i = 1,2, \ldots, p$.}
\item {\bf Variational problem}: For a given $\epsilon \in (0, \epsilon_{\rm max})$, determine an admissible optimal control $u^{\epsilon}(t)$, $t \in [0, T]$, that minimizes the following functional
 \begin{align*}
 J^{\epsilon}[u] = \Phi\bigl(\theta^{\epsilon}(T), \mathcal{Z}^{(2)}\bigr), \quad \text{s.t. ~~~ Equation}~\eqref{Eq1.2},
\end{align*}
 where $\Phi\bigl(\theta, \mathcal{Z}^{(2)}\bigr)$ is a scalar function that depends on the model validating dataset $\mathcal{Z}^{(2)}$. Note that such a variational problem together with the above weakly-controlled gradient system provides a mathematical apparatus how to improve generalization in such a class of learning problems.\footnote{In this paper, we consider the following function
 \begin{align*}
\Phi\bigl(\theta^{\epsilon}(T), \mathcal{Z}^{(2)}\bigr) = \big(1/m_2\big) \sum\nolimits_{i=1}^{m_2} {\ell} \bigl(h_{\theta^{\epsilon}(T)}(x_i^{(2)}), y_i^{(2)} \bigr), \quad \text{w.r.t. the dataset} ~ \mathcal{Z}^{(2)},
\end{align*}
as a cost functional $J^{\epsilon}[u]$ that serves as a measure for evaluating the quality of the estimated optimal parameter $\theta^{\ast} = \theta^{\epsilon}(T)$, i.e., when $\theta^{\epsilon}(t)$ is evaluated at a certain fixed time $T$.}
 \item {\bf General assumptions:} Throughout the paper, we assume the following conditions: (i) the set $U$ is compact in $\mathbb{R}$ and the final time $T$ is fixed, (ii) the function $\Phi\bigl(\theta, \mathcal{Z}^{(2)}\bigr)$ is twice continuously differentiable w.r.t. the parameter $\theta$, and (iii) for any $\epsilon \in (0, \epsilon_{\rm max})$ and all admissible bounded controls $u(t)$ from $U$, the solution ${\theta}^{\epsilon}(t)$ for $t \in [0, T]$ which corresponds to the weakly-controlled gradient system of Equation~\eqref{Eq1.2} starting from an initial condition $\theta^{\epsilon}(0)=\theta_0$, exists and bounded.\footnote{For such a variational problem, these assumptions are sufficient for the existence of a nonempty compact reachable set $\mathcal{R}(\theta_0) \subset \Theta$, for some admissible controls on $[0,T]$ that belongs to $U$, starting from an initial point $\theta^{\epsilon}(0)=\theta_0$ (e.g., see \cite{r1} for related discussions on the Filippov's theorem providing a sufficient condition for compactness of the reachable set).}.
 \end{enumerate}

In what follows, we assume that there exists an admissible optimal control $u^{\epsilon}(t)$ for all $\epsilon \in (0, \epsilon_{\rm max})$. Then, the necessary optimality conditions for the optimal control problem with weakly-controlled gradient system satisfy the following Euler-Lagrange critical point equations
\begin{align}
 \dot{\theta}^{\epsilon}(t) &= \frac{\partial H^{\epsilon}\bigl(\theta^{\epsilon}(t), p^{\epsilon}(t), u^{\epsilon}(t)\bigr)}{\partial p}, \notag \\
                                       &= - \nabla J_{0} \bigl(\theta^{\epsilon}(t), \mathcal{Z}^{(1)}\bigr) + \epsilon u^{\epsilon}(t) B\bigl(\theta^{\epsilon}(t), \tilde{\mathcal{Z}}^{(1)}\bigr),  \quad \theta^{\epsilon}(0)=\theta_0, \label{Eq1.3} \\
   \dot{p}^{\epsilon}(t) &= -\frac{\partial H^{\epsilon}\bigl(\theta^{\epsilon}(t), p^{\epsilon}(t), u^{\epsilon}(t)\bigr)}{\partial \theta}. \notag \\
                                       &= \nabla^2 J_{0} \bigl(\theta^{\epsilon}(t), \mathcal{Z}^{(1)}\bigr) p^{\epsilon}(t) - \epsilon u^{\epsilon}(t) \bigl(\nabla B\bigl(\theta^{\epsilon}(t), \tilde{\mathcal{Z}^{(1)}}\bigr) \bigr)^T p^{\epsilon}(t), \notag\\
                                       & \quad \quad \quad \quad \quad  \quad p^{\epsilon}(T) = - \nabla \Phi\bigl(\theta^{\epsilon}(T), \mathcal{Z}^{(2)}\bigr), \label{Eq1.4}\\
                                       u^{\epsilon}(t) &= \argmax H^{\epsilon}\bigl(\theta^{\epsilon}(t), p^{\epsilon}(t), u(t)\bigr), ~~ u(t) \in U ~~ \text{on} ~~ t \in [0, T], \label{Eq1.5}
\end{align}
where the Hamiltonian function $H^{\epsilon}$ is given by
\begin{align}
 H^{\epsilon} \bigl(\theta, p, u \bigr) = \bigl \langle p,\, - \nabla J_{0} \bigl(\theta, \mathcal{Z}^{(1)}\bigr) + \epsilon u B\bigl(\theta, \tilde{\mathcal{Z}}^{(1)}\bigr) \bigr \rangle \label{Eq1.6}
\end{align}
and such optimality conditions are the direct consequence of the {\it Pontryagin's maximum principle} (e.g., see \cite{r2} for additional discussions on the first-order necessary optimality conditions; see also \cite{r3} for related discussions in the context of learning). 

In following section, using the perturbation theory (e.g., see \cite{r4}, \cite{r5} or \cite{r6} for additional discussions), we provide approximation solutions for the variational problem related to the weakly-controlled gradient system with a small parameter. In particular, we provide results that will allow us to solve independently a sequence of optimization problems, i.e., a set of decomposed optimization problems, so as to aggregate the corresponding approximate optimal solutions that are reasonably sufficient for improving generalization in such a class of learning problems. Moreover, we also provide an estimate for its rate of convergence.

\section{Main results} \label{S2}
Assume that the solutions $\theta^{\epsilon}(t)$,  $p^{\epsilon}(t)$ and $u^{\epsilon}(t)$ corresponding to the optimality conditions in Equations~\eqref{Eq1.3}-\eqref{Eq1.5} can be expressed as a series in the small parameter $\epsilon$ as follows
\begin{equation}
\left.\begin{matrix} 
\theta^{\epsilon}(t) &= \theta^{0}(t) + \epsilon \theta^{1}(t) + \mathcal{O}\bigl(\epsilon^2\bigr)\\
p^{\epsilon}(t) &= p^{0}(t) + \epsilon p^{1}(t) + \mathcal{O}\bigl(\epsilon^2\bigr)\\
u^{\epsilon}(t) &= u^{0}(t) + \epsilon u^{1}(t) + \mathcal{O}\bigl(\epsilon^2\bigr)
\end{matrix}\right\}  \label{Eq2.1}
\end{equation}  
Then, the following proposition (whose proof is given in the Appendix section) characterizes the approximate solutions obtained by keeping only the first two terms (i.e., the zeroth and the first-order solutions).
\begin{proposition} \label{P1}
Let the optimal solutions $\theta^{\epsilon}(t)$, $p^{\epsilon}(t)$ and $u^{\epsilon}(t)$ be expressed as a series in the small parameter $\epsilon$ (i.e., as in Equation~\eqref{Eq2.1} above). Then, the optimal value for the variational problem, i.e.,
\begin{align*}
J^{\epsilon}[u] &= \Phi\bigl(\theta^{\epsilon}(T), \mathcal{Z}^{(2)}\bigr) \quad \to \quad \min_{u(t) \in U ~ \text{on} ~ t \in [0, T]}, \quad \text{s.t. ~~~ Equation}~\eqref{Eq1.2},
\end{align*}
satisfies the following condition
\begin{align}
J^{\epsilon}[u^{\epsilon}] &= \Phi\bigl(\theta^{0}(T), \mathcal{Z}^{(2)}\bigr) + \epsilon \bigl \langle \nabla \Phi\bigl(\theta^{0}(T), \mathcal{Z}^{(2)}\bigr), \, \theta^{1}(T) \bigr\rangle + \mathcal{O}\bigl(\epsilon^2\bigr), \notag\\
                                      &= \Phi\bigl(\theta^{0}(T), \mathcal{Z}^{(2)}\bigr) - \epsilon \bigl \langle p^{0}(T), \, \theta^{1}(T) \bigr\rangle + \mathcal{O}\bigl(\epsilon^2\bigr), \label{Eq2.2}
\end{align}
where the zeroth-order solutions $\theta^{0}(t)$,  $p^{0}(t)$ and $u^{0}(t)$ satisfy the following critical conditions
\begin{align}
 \dot{\theta}^{0}(t) &= - \nabla J_{0} \bigl(\theta^{0}(t), \mathcal{Z}^{(1)}\bigr),  \quad \theta^{0}(0)=\theta_0, \label{Eq2.3} \\
   \dot{p}^{0}(t) &= \nabla^2 J_{0} \bigl(\theta^{0}(t), \mathcal{Z}^{(1)}\bigr), \quad p^{0}(T) = - \nabla \Phi\bigl(\theta^{0}(T), \mathcal{Z}^{(2)}\bigr), \label{Eq2.4} \\
           u^{0}(t) &= \argmax \bigl \langle p^{0}(t), u(t) B\bigl(\theta^{0}(t), \tilde{\mathcal{Z}}^{(1)}\bigr)\bigr \rangle, ~~ u(t) \in U ~~\text{on} ~~ t \in [0, T], \label{Eq2.5}
\end{align}
and the first-order solution $\theta^{1}(t)$ satisfies the following system dynamics equation
\begin{align}
 \dot{\theta}^{1}(t) = - \nabla^2 J_{0} \bigl(\theta^{0}(t), \mathcal{Z}^{(1)}\bigr) \theta^{1}(t) + u^{0}(t) B\bigl(\theta^{0}(t), \tilde{\mathcal{Z}}^{(1)}\bigr), ~~ \theta^{1}(0)=0, ~~ t \in [0, T]. \label{Eq2.6}
\end{align}
Moreover, the estimated optimal parameter $\theta^{\ast}$ is given by
\begin{align}
\theta^{\ast} &= \theta^{\epsilon}(T), \notag\\
                    &= \theta^{0}(T) + \epsilon \theta^{1}(T) + \mathcal{O}\bigl(\epsilon^2\bigr). \label{Eq2.7}
\end{align}
\end{proposition}

Here, we remark that $\theta^{1}(t)$ depends on the the zeroth-order solutions $\theta^{0}(t)$, $p^{0}(t)$ and $u^{0}(t)$, for $t \in [0,T]$, as well as on the datasets $\mathcal{Z}^{(1)}$, $\mathcal{Z}^{(2)}$ and $\tilde{\mathcal{Z}}^{(1)}$ (where the latter is generated by adding random noise with small variance $\sigma^2$ in the model training dataset $\mathcal{Z}^{(1)}$). Note that the control $u^{0}(t)$, for $t \in [0, T]$, may not be unique, but it is measurable and belongs to $U$ for all $t \in [0, T]$.\footnote{Note that if we are looking for an optimal control $u(t) \in U$, for $t \in [0, T]$, from a class of bounded functions $U=\big\{u \in \mathbb{R} \colon -1 \le u \le +1 \big\}$. Then, the solution satisfying the extremum condition of Equation~\eqref{Eq2.5} takes values between $-1$ or $+1$, i.e., a {\it bang-bang control} with an interpretation of switching which is consistent with nonsmooth calculus of variations (see \cite{r7} or \cite{r8} for additional discussions).} Moreover, the improvements in the model training loss and that of the model validating loss are given by 
\begin{align}
 &J_{0} \bigl(\theta^{\epsilon}(T), \mathcal{Z}^{(1)}\bigr) - J_{0} \bigl(\theta^{0}(T), \mathcal{Z}^{(1)}\bigr) \notag \\
    &  \quad \quad = \frac{1}{m_1} \sum\nolimits_{i=1}^{m_1} \bigl({\ell} \bigl(h_{\theta^{0}(T)+ \epsilon \theta^{1}(T)}(x_i^{(1)}), y_i^{(1)} \bigr) - {\ell} \bigl(h_{\theta^{0}(T)}(x_i^{(1)}), y_i^{(1)} \bigr)\bigr)\label{Eq2.8}
\end{align}
and
\begin{align}
 &J_{0} \bigl(\theta^{\epsilon}(T), \mathcal{Z}^{(2)}\bigr)  - J_{0} \bigl(\theta^{0}(T), \mathcal{Z}^{(2)}\bigr)  \notag \\
       &\quad \quad  = \frac{1}{m_2} \sum\nolimits_{i=1}^{m_2} \bigl({\ell} \bigl(h_{\theta^{0}(T)+ \epsilon \theta^{1}(T)}(x_i^{(2)}), y_i^{(2)} \bigr) - {\ell} \bigl(h_{\theta^{0}(T)}(x_i^{(2)}), y_i^{(2)} \bigr)\bigr), \label{Eq2.9}
\end{align}
respectively.

Finally, the following proposition provides an estimate for the convergence rate (see the Appendix section for its proof).
\begin{proposition} \label{P2}
 Suppose that Proposition~\ref{P1} holds, then we have
 \begin{align}
 J^{\epsilon}[u^{\epsilon}] - J^{\epsilon}[u^{0}] = \mathcal{O}\bigl(\epsilon^2\bigr) \label{Eq2.10}
\end{align}
 for the estimate of convergence rate.
 \end{proposition}
 
 \section*{Algorithm} Here, we present a generic algorithm for solving such approximate optimal solutions - based on the zeroth and the first-order solutions - for the optimal control problem in Equations~\eqref{Eq1.1} and \eqref{Eq1.2}.

{\rm \footnotesize

{\bf ALGORITHM:}
\begin{itemize}
\item[{\bf 1.}] {\bf The zeroth-order solutions:} Solve the forward and backward-equations w.r.t. the system dynamics of Equations~\eqref{Eq2.3} and \eqref{Eq2.4}, i.e.,
\begin{align*}
 \dot{\theta}^{0}(t) &= - \nabla J_{0} \bigl(\theta^{0}(t), \mathcal{Z}^{(1)}\bigr),  \quad \theta^{0}(0)=\theta_0, \\
   \dot{p}^{0}(t) &= \nabla^2 J_{0} \bigl(\theta^{0}(t), \mathcal{Z}^{(1)}\bigr), \quad p^{0}(T) = - \nabla \Phi\bigl(\theta^{0}(T), \mathcal{Z}^{(2)}\bigr).
\end{align*}
\item[{\bf 2.}] Compute the admissible control $u^{(0)}(t) \in U$ for $t \in [0, T]$ using Equation~\eqref{Eq2.5}, i.e.,
\begin{align*}
     u^{0}(t) = \argmax \bigl \langle p^{0}(t), u(t) B\bigl(\theta^{0}(t), \tilde{\mathcal{Z}}^{(1)}\bigr)\bigr \rangle, ~~ u(t) \in U ~~\text{on} ~~ t \in [0, T].
\end{align*}
\item[{\bf 3.}] {\bf The first-order solutions:} Using $\theta^{0}(t)$, $p^{0}(t)$ and $u^{0}(t)$, for $t \in [0, T]$, solve the forward system equation corresponding to Equation~\eqref{Eq2.6}, i.e.,
 \begin{align*}
 \dot{\theta}^{1}(t) = - \nabla^2 J_{0} \bigl(\theta^{0}(t), \mathcal{Z}^{(1)}\bigr) \theta^{1}(t) + u^{0}(t) B\bigl(\theta^{0}(t), \tilde{\mathcal{Z}}^{(1)}\bigr), ~~ \theta^{1}(0)=0, ~~ t \in [0, T].
\end{align*}
 \item[{\bf 4.}] {\bf Output:} Return the estimated optimal parameter value $\theta^{\ast} = \theta^{0}(T) + \epsilon \theta^{1}(T)$.
\end{itemize}}

Here, it is worth remarking that the above generic algorithm allows us to determine the zeroth and the first-order solutions, corresponding to a sequence of decomposed optimization problems, and their aggregations are reasonable approximate optimal solutions to the optimal control problem in Equations~\eqref{Eq1.1} and \eqref{Eq1.2}.

\section{Numerical results and discussions}  \label{S3}
In this section, we presented some numerical results for a simple polynomial interpolation problem modeling the thermophysical properties of saturated water (in liquid state), where the dataset (which is given in Table~\ref{Tb3} of the Appendix section) is taken from (see \cite[p.~1003]{r9}). Here, we only considered the problem of point estimation for modeling of (i) the density $\rho$, (ii) the specific heat $c_{p}$, and (iii) the thermal conductivity $k$, where the mathematical model relating each of these thermophysical properties as a function of the temperature ${\rm T}$ (in Kelvin ${\rm [K]}$, in the range $273.15 \,{\rm K} \le {\rm T} \le 373.15\, {\rm K}$), is assumed to obey a second-order polynomial function of the form $h_{\theta}({\rm T}) = \theta_1 + \theta_2 {\rm T} + \theta_3 {\rm T}^2$, with $\theta = \bigl(\theta_1, \theta_2, \theta_3\bigr)$ is the model parameter.\footnote{The notation for the temperature ${\rm T}$ is different from that of the fixed final time $T$ which is part of the problem statement.}

\begin{figure}[hbt]
\begin{center}
 \includegraphics[scale=0.125]{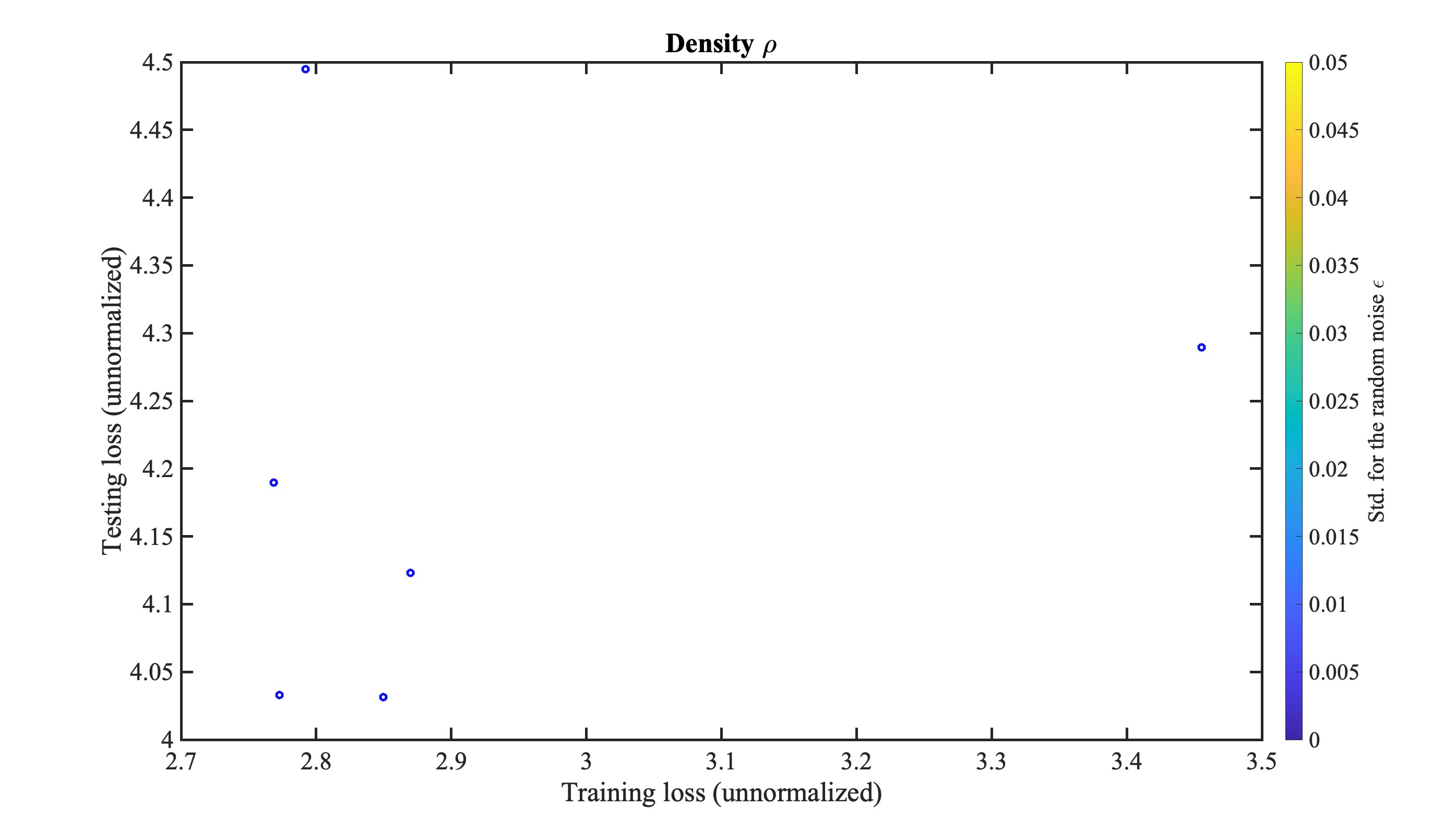}
 \caption{Plot for model training loss versus model testing loss for different levels of random noise distortions (Density $\rho$).} \label{Fig1}
\end{center}
\end{figure}

\begin{figure}[hbt]
\begin{center}
 \includegraphics[scale=0.125]{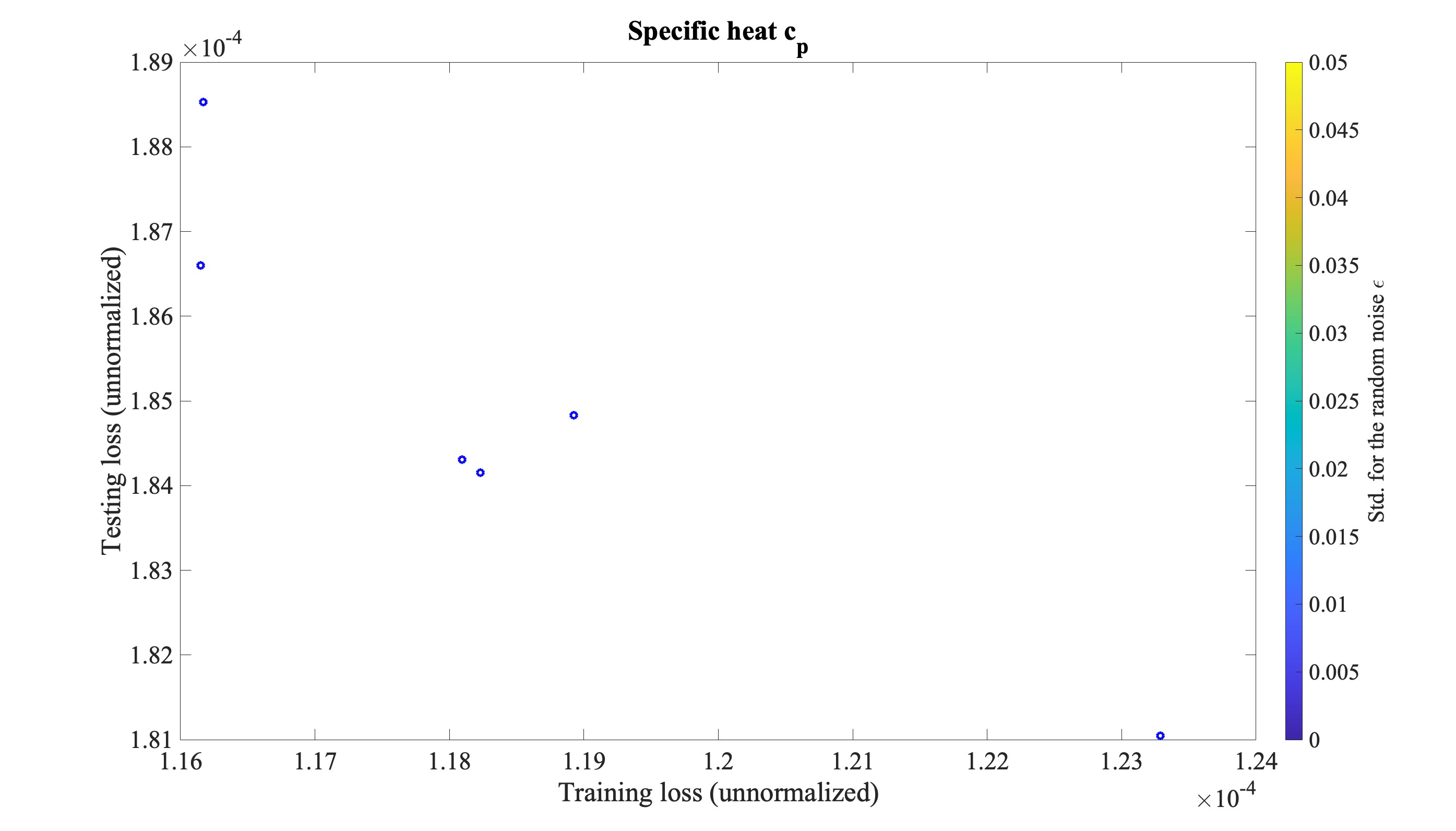}
  \caption{Plot for model training loss versus model testing loss for different levels of random noise distortions (Specific heat $c_p$).} \label{Fig2}
\end{center}
\end{figure}

\begin{figure}[hbt]
\begin{center}
   \includegraphics[scale=0.125]{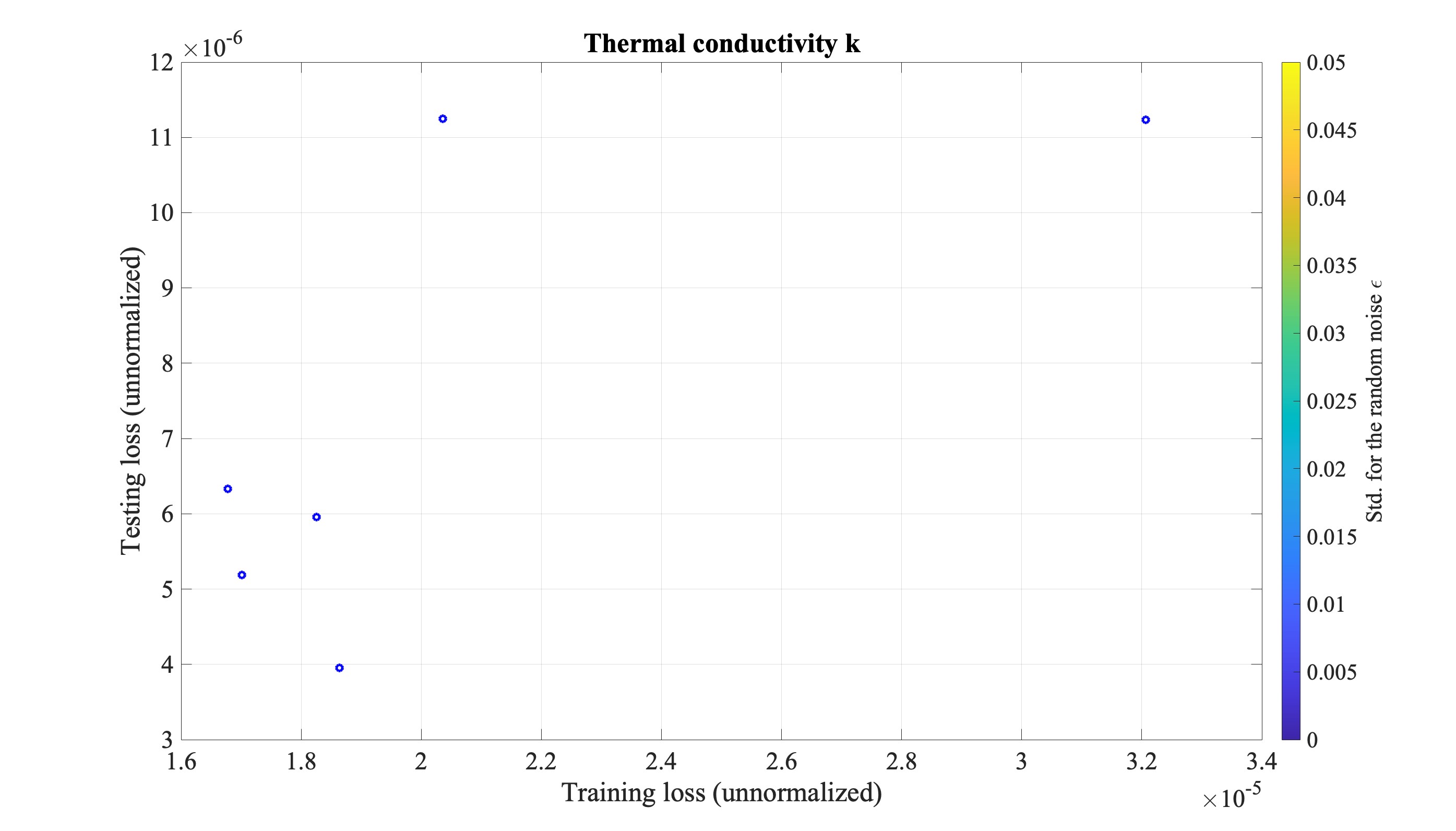}
  \caption{Plot for model training loss versus model testing loss for different levels of random noise distortions (Thermal conductivity $k$).} \label{Fig3}
\end{center}
\end{figure}

For the numerical simulation result, we first partitioned the original dataset $\mathcal{Z}^{(0)}$ of size $m_0=22$ into two data subsets, i.e., model training dataset $\mathcal{Z}^{(1)}$ and model validating dataset $\mathcal{Z}^{(2)}$, of sizes $m_1 = 18$ and $m_2 = 6$, respectively. Moreover, for the dataset $\tilde{\mathcal{Z}}^{(1)}$, which is associated with the nonlinear term $B$ in the system dynamics, we added random noises with distortion levels of $1\%$ and $5\%$ corresponding to the sample variance of the original dataset and with a value of small parameter $\epsilon = 0.001$. Tables~\ref{Tb1} and \ref{Tb2} show the estimated optimal parameters $\theta_i^{\ast} = \theta_i^{0}(T) + \epsilon \theta_i^{1}(T)$, for $i=1,2,3$, with fixed final time $T$, and the sample standard deviation for the residual errors $\varepsilon^{\rm res}$. Note that we computed the residual errors w.r.t. the original dataset $\mathcal{Z}^{(0)}$ based on $\varepsilon_i^{\rm res} = y_i - h_{\theta^{\ast}}({\rm T}_i)$, for $i=1,2, \ldots, 22$. Finally, Figures~\ref{Fig1}-\ref{Fig3} show the model training loss versus that of the model testing loss for different levels of random noise distortions.

\newpage

\newpage

\begin{table}[hbt]
\begin{center}
{\footnotesize
\caption{The estimated optimal parameters with a level of $1\%$ random noise distortion} \label{Tb1}
  \begin{tabular}{c | c c c | c }\hline
   & & & & {\rm Sample} \\
  $h_{\theta}({\rm T})$  & $\theta_1^{\ast}$ & $\theta_2^{\ast}$ & $\theta_3^{\ast}$ & {\rm standard deviation for}  \\
    & & & & {\rm the residual errors $\varepsilon^{\rm res}$}  \\ \hline
        $\rho$ & $763.1823$ & $1.8221$ & $-3.4862 \times 10^{-3}$ & $0.5693$  \\
        $c_p$ & $5.5944$ & $-8.8978 \times 10^{-3}$  & $1.3982 \times 10^{-5}$ & $0.0038$  \\
        $k$ & $-0.4338$ & $0.0056$ & $-6.9164 \times 10^{-6}$ & $0.0010$ \\\hline
   \end{tabular}}
  \end{center}
 \end{table}

\begin{table}[hbt]
\begin{center}
{\footnotesize
\caption{The estimated optimal parameters with a level of $5\%$ random noise distortion} \label{Tb2}
  \begin{tabular}{c | c c c | c }\hline
   & & & & {\rm Sample}  \\
  $h_{\theta}({\rm T})$  & $\theta_1^{\ast}$ & $\theta_2^{\ast}$ & $\theta_3^{\ast}$ & {\rm standard deviation for}  \\
    & & & & {\rm the residual errors $\varepsilon^{\rm res}$}  \\ \hline
        $\rho$ & $759.7205$ & $1.8512$ & $-3.5432 \times 10^{-3}$ & $0.6221$  \\
        $c_p$ & $5.6124$ & $-9.0082 \times 10^{-3}$  & $1.4148 \times 10^{-5}$ & $0.0038$  \\
        $k$ & $-0.4769$ & $0.0058$ & $-7.3308 \times 10^{-6}$ & $0.0011$ \\\hline
   \end{tabular}}
  \end{center}
 \end{table}

\section*{A.~Appendix}
{\it Proof of Proposition~\ref{P1}.}
Clearly, if the expansions of Equation~\eqref{Eq2.1} in the small parameter $\epsilon$ converge and satisfy the Euler-Lagrange critical point equations of Equations~\eqref{Eq1.3}-\eqref{Eq1.5}. Then, with direct substitutions, we have the following relations
\begin{align*}
 \dot{\theta}^{0}(t) + \epsilon \dot{\theta}^{1}(t) + \mathcal{O}\bigl(\epsilon^2\bigr) &= - \nabla J_{0} \bigl(\theta^{0}(t), \mathcal{Z}^{(1)}\bigr) - \epsilon \nabla^2 J_{0} \bigl(\theta^{0}(t), \mathcal{Z}^{(1)}\bigr) \theta^{1}(t) \\
 & \quad \quad \quad + \epsilon u^{0}(t) B\bigl(\theta^{0}(t), \tilde{\mathcal{Z}}^{(1)}\bigr) + \mathcal{O}\bigl(\epsilon^2\bigr), \quad \theta^{0}(0)=\theta_0
 \end{align*}
 and
\begin{align*}
  \dot{p}^{0}(t) +\epsilon \dot{p}^{1}(t) + \mathcal{O}\bigl(\epsilon^2\bigr) &= \nabla^2 J_{0} \bigl(\theta^{0}(t), \mathcal{Z}^{(1)}\bigr) p^{0}(t) + \epsilon \nabla^2 J_{0} \bigl(\theta^{0}(t), \mathcal{Z}^{(1)}\bigr) p^{1}(t), \\
  & \quad \quad - \epsilon u^{0}(t) \bigl(\nabla B\bigl(\theta^{0}(t), \tilde{\mathcal{Z}}\bigr)^T p^{1}(t) + \mathcal{O}\bigl(\epsilon^2\bigr), \\
  & \quad \quad \quad \quad p^{0}(T) = - \nabla \Phi\bigl(\theta^{0}(T), \mathcal{Z}^{(2)}\bigr).
\end{align*}
where the admissible control $u^{0}(t) \in U$, for $ t \in [0, T]$, is an extremum solution of 
\begin{align*}
u^{0}(t) &= \argmax \bigl \langle p^{0}(t), u(t) B\bigl(\theta^{0}(t), \tilde{\mathcal{Z}}^{(1)}\bigr)\bigr \rangle, ~~ u(t) \in U ~~\text{on} ~~ t \in [0, T].
\end{align*}
Note that if we equate equal powers of $\epsilon$ and only retaining the zeroth and first-order solutions. Then, we will arrive to Equations~\eqref{Eq2.3} and \eqref{Eq2.4}. Moreover, the optimal estimated parameter $\theta^{\ast} = \theta^{\epsilon}(T)$ can be recovered from the following relation
\begin{align*}
\theta^{\epsilon}(T) = \theta^{0}(T) + \epsilon \theta^{1}(T) + \mathcal{O}\bigl(\epsilon^2\bigr). 
\end{align*}
This completes the proof of Proposition~\ref{P1}.\qed

{\it Proof of Proposition~\ref{P2}.}
Note that, for any admissible control $u(t) \in U$, for $t \in [0, T]$, we have the following 
\begin{align*}
 J^{\epsilon}[u] &= \Phi\bigl(\theta^{0}(T), \mathcal{Z}^{(2)}\bigr) + \epsilon \bigl \langle \nabla \Phi\bigl(\theta^{0}(T), \mathcal{Z}^{(2)}\bigr), \, \theta^{1}(T) \bigr\rangle + \mathcal{O}\bigl(\epsilon^2\bigr), \\
                &= \Phi\bigl(\theta^{0}(T), \mathcal{Z}^{(2)}\bigr) - \epsilon \bigl \langle p^{0}(T), \, \theta_{u}^{1}(T) \bigr\rangle + \mathcal{O}\bigl(\epsilon^2\bigr),
\end{align*}
which is uniform w.r.t. $u(t)$ and $p^{0}(T)$, where $\theta_{u}^{1}(T)$ is the solution of
\begin{align*}
 \dot{\theta}_{u}^{1}(t) = - \nabla^2 J_{0} \bigl(\theta^{0}(t), \mathcal{Z}^{(1)}\bigr) \theta_{u}^{1}(t) + u(t) B\bigl(\theta^{0}(t), \tilde{\mathcal{Z}}^{(1)}\bigr), ~~\text{with} ~~ \theta_{u}^{1}(0)=0
\end{align*}
evaluated at $t = T$, and $\theta^{0}(T)$ is the solution of Equation~\eqref{Eq2.3} evaluated at $t = T$, while $p^{0}(T) = - \nabla \Phi\bigl(\theta^{0}(T), \mathcal{Z}^{(2)}\bigr)$. 

As a result, we have the following relation, w.r.t. $u^{\epsilon}(t)$ and $u^{0}(t)$ for $t \in [0, T]$,
\begin{align*}
 J^{\epsilon}[u^{\epsilon}] - J^{\epsilon}[u^{0}] &= \epsilon \bigl \langle p^{0}(T), \, \theta_{u^{0}}^{1}(T) - \theta_{u^{\epsilon}}^{1}(T) \bigr \rangle + \mathcal{O}\bigl(\epsilon^2\bigr),\\
   &= \epsilon \int_0^T \bigl \langle \dot{p}^{0}(t), \, \theta_{u^{0}}^{1}(t) - \theta_{u^{\epsilon}}^{1}(t) \bigr \rangle dt \\
   & \quad \quad \quad + \epsilon \int_0^T \bigl \langle p^{0}(t), \, \dot{\theta}_{u^{0}}^{1}(t) - \dot{\theta}_{u^{\epsilon}}^{1}(t) \bigr\rangle dt + \mathcal{O}\bigl(\epsilon^2\bigr).
\end{align*}
Moreover, if we substitute the following two system equations
\begin{align*}
 \dot{p}^{0}(t) &= \nabla^2 J_{0} \bigl(\theta^{0}(t), \mathcal{Z}^{(1)}\bigr)
\end{align*}
and
\begin{align*}
 \dot{\theta}_{u^{0}}^{1}(t) - \dot{\theta}_{u^{\epsilon}}^{1}(t) &= - \nabla^2 J_{0} \bigl(\theta^{0}(t), \mathcal{Z}^{(1)}\bigr) \bigl(\theta_{u^{0}}^{1}(t) - \theta_{u^{\epsilon}}^{1}(t)\bigr) \\
 & \quad \quad \quad \quad + \bigl(u^{0}(t) - u^{\epsilon}(t) \bigr)B\bigl(\theta^{0}(t), \tilde{\mathcal{Z}}^{(1)}\bigr)
\end{align*}
in the equation of $J^{\epsilon}[u^{\epsilon}] - J^{\epsilon}[u^{0}]$ above. Then, we have the following relation
\begin{align*}
 J^{\epsilon}[u^{\epsilon}] - J^{\epsilon}[u^{0}] = \epsilon \bigl(u^{0}(t) - u^{\epsilon}(t) \bigr)B\bigl(\theta^{0}(t), \tilde{\mathcal{Z}}^{(1)}\bigr) + \mathcal{O}\bigl(\epsilon^2\bigr).
 \end{align*}
Furthermore, noting that $\bigl(u^{0}(t) - u^{\epsilon}(t) \bigr) = \mathcal{O}\bigl(\epsilon\bigr)$ (cf. Equation~\eqref{Eq2.1}), then we will arrive to 
\begin{align*}
 J^{\epsilon}[u^{\epsilon}] - J^{\epsilon}[u^{0}] = \mathcal{O}\bigl(\epsilon^2\bigr).
 \end{align*}
 for the estimate of convergence rate. This completes the proof of Proposition~\ref{P2}.\qed

\section*{B.~Appendix}~\\

\begin{table}[h]
\begin{center}
{\footnotesize
\caption{Thermophysical properties of saturated water (Liquid).} \label{Tb3}
  \begin{tabular}{c c c c c c}\hline
           ${\rm T}$  & $\rho$ & $c_p$  & $\mu \times 10^{6}$ & $k$  & $\beta \times 10^{6}$ \\
  $({\rm K})$  & $({\rm kg/m^3})$ & $({\rm kJ/kg \cdot K})$ & $({\rm N \cdot s/m^2})$ & $({\rm W/m\cdot K })$ & $({\rm K^{-1}})$  \\ \hline
        273.15 & 1000 & 4.217 & 1750 & 0.569 & -68.05\\
        273 & 1000 & 4.211 & 1652 & 0.574 & -32.74\\
        280 & 1000 & 4.198 & 1422 & 0.582 & 46.04\\
        285 & 1000 & 4.189 & 1225 & 0.590 & 114.1\\
        290 &   999 & 4.184 & 1080 & 0.598 & 174.0\\
        295 &   998 & 4.181 &   959 & 0.606 & 227.5\\
        300 &   997 & 4.179 &   855 & 0.613 & 276.1\\
        305 &   995 & 4.178 &   769 & 0.620 & 320.6\\
        310 &   993 & 4.178 &   695 & 0.628 & 361.9\\
        315 &   991 & 4.179 &   631 & 0.634 & 400.4\\
        320 &   989 & 4.180 &   577 & 0.640 & 436.7\\
        325 &   987 & 4.182 &   528 & 0.645 & 471.2\\
        330 &   984 & 4.184 &   489 & 0.650 & 504.0\\
        335 &   982 & 4.186 &   453 & 0.656 & 535.5\\
        340 &   979 & 4.188 &   420 & 0.660 & 566.0\\
        345 &   977 & 4.191 &   389 & 0.664 & 595.4\\
        350 &   974 & 4.195 &   365 & 0.668 & 624.2\\
        355 &   971 & 4.199 &   343 & 0.671 & 652.3\\
        360 &   967 & 4.203 &   324 & 0.674 & 679.9\\
        365 &   963 & 4.209 &   306 & 0.677 & 707.1\\
        370 &   961 & 4.214 &   289 & 0.679 & 728.7\\
        373.15 & 958  & 4.217 &   279 & 0.680 & 750.1\\
        400 &   937 & 4.256 &   217 & 0.688 & 896\\
        450 &   890 & 4.40   &  152 & 0.678 & \\
        500 &   831 & 4.66   &  118 & 0.642 & \\
        550 &   756 & 5.24   &    97 & 0.580 & \\
        600 &   649 & 7.00   &    81 & 0.497 & \\
        647.3 &  315 & 0      &    45 & 0.238 & \\ \hline
  \end{tabular}}
  \end{center}
 \end{table}

\end{document}